\documentstyle{article}
\def\Bbb{\bf}

\def\be{\begin{equation}}
\def\ee{\end{equation}}
\def\bea{\begin{eqnarray*}}
\def\eea{\end{eqnarray*}}
\newtheorem{main}{Theorem}

\newtheorem{thm}{Theorem}
\newtheorem{prop}{Proposition}
\newtheorem{cor}{Corollary}

\newenvironment{proof}{\medskip \noindent
{\bf Proof.}}{\hfill \rule{.5em}{1em} \medskip}
\begin{document}
\sloppy
\title{Yamabe Constants and the Perturbed \\Seiberg-Witten Equations}

\author{Claude LeBrun\thanks{Supported 
in part by  NSF grant DMS-9505744.} 
\\ 
SUNY Stony
 Brook 
  }

\date{}
\maketitle

\begin{abstract} 
Among all conformal classes of Riemannian 
metrics on    ${\Bbb CP}_2$,
that of the Fubini-Study metric is shown
to have the largest Yamabe constant. The proof, which involves  
perturbations of the  Seiberg-Witten equations, 
also yields  new results on the total scalar curvature
of almost-K\"ahler 4-manifolds. 
 \end{abstract}
 
 \section{Introduction}
One may define an interesting and natural diffeomorphism invariant
of a compact smooth $n$-manifold $M$  as a  minimax 
of   the total  scalar curvature
over all unit volume Riemannian metrics on $M$.
To be precise, for each    Riemannian metric $g$ on $M$, one
 may first 
define an invariant of the corresponding  conformal class
$[g]=\{ vg ~|~v: M\to {\Bbb R}^+\}$ by setting 
$$Y_{[g]} = \inf_{g\in [g] } \frac{\int_M 
s_g~d\mu_g}{\left(\int_M 
d\mu_g\right)^{\frac{n-2}{n}}},$$
where $s_g$ and $d\mu_g$ respectively
 denote the scalar curvature and volume $n$-form of
$g$. The number $Y_{[g]}$ is called the {\em Yamabe constant}
of the conformal class of metrics, and a remarkable theorem
\cite{aubin,lp,rick} 
of Yamabe, Trudinger, Aubin, and Schoen states that any conformal class 
$[g]$ contains a metric (necessarily of constant scalar curvature) 
which achieves this infimum.  
 We may now define an invariant
$$Y(M) = \sup_{[g]} Y_{[g]} = \sup_{[g]}\inf_{g\in [g] } \frac{\int_M
s_g~d\mu_g}{\left(\int_M 
d\mu_g\right)^{\frac{n-2}{n}}}$$
which we will call the {\em Yamabe invariant} of $M$.

For 2-manifolds, this invariant is easy to compute;
 the Gauss-Bonnet theorem tells us
that the Yamabe invariant of  a compact surface of
genus $\bf g$ is just $8\pi(1-{\bf g})$. 
But relatively few precise computations of Yamabe invariants have
proved tractable  for  manifolds of higher dimension.
One important general fact, due to Aubin, is that 
for any compact $n$-manifold one has 
$Y(M)\leq Y(S^n)= n(n-1) V^{2/n}_n$, where $V_n$ is the 
volume of the unit $n$-sphere.  Arguments  
involving the Dirac operator \cite{GL,bes} have also be used 
for some time to show that $Y=0$
for the $n$-torus and certain other manifolds admitting Ricci-flat
metrics. It was nonetheless a revolutionary development when 
Witten \cite{witten} observed that a non-linear Dirac sytem 
called the  {\em Seiberg-Witten equations} can allow one
to distinguish different smooth structures on
a  topological 4-manifold by virtue of the fact 
that the corresponding Yamabe invariants are distinct.  
Pushing this insight to its logical conclusion
\cite{leb3} results in precise computations
of  the Yamabe invariant for an infinite class  of 4-manifolds
with $Y(M) < 0$.  

The usual Seiberg-Witten invariant vanishes, however, when
$Y(M) > 0$, so it might seem hopeless to try to 
use such techniques to show that certain  4-manifolds 
have $0< Y(M) < Y(S^4) = 8\sqrt{6}\pi$. For
4-manifolds with $b_+=1$, however, it turns out that 
a perturbed version of the equations \cite{taubes,leb1,liu}
gives rise to an invariant which is often non-trivial
even in the presence of metrics of positive scalar
curvature. Extracting information about the scalar curvature
becomes more complicated, however, since the perturbation 
dominates the metric terms in
the corresponding Weitzenb\"ock formula. Nonetheless, if 
the perturbed Seiberg-Witten equations on $(M^4,g)$ have
a solution when the
 perturbation term is 
any large multiple  of the self-dual harmonic form
$\omega$, then one obtains (\S \ref{key}) the basic estimate 
$$\int s \frac{|\omega |}{\sqrt{2}} d\mu 
\leq 4\pi c_1\cdot [\omega
],$$
with  equality  iff
$g$ is  K\"ahler. 
The presence of $\omega$ in this estimate 
exactly compensates for the fact that the 
Yamabe functional is  neither bounded above nor bounded  below.
One  eventually estimates the Yamabe constants
$Y_{[g]}$ by (\S \ref{rep}) representing each  conformal class 
by metrics for which the size of $\omega$ is
nearly constant.  
The main application is the following:

\begin{main}
The Yamabe invariant of the complex projective plane is given by
$$Y({\Bbb CP}_2) = 12\sqrt{2}\pi .$$
Moreover, a conformal class $[g]$ on ${\Bbb CP}_2$ satisfies
$Y_{[g]} = Y({\Bbb CP}_2)$ iff there is a diffeomorphism 
$\Phi : {\Bbb CP}_2\to {\Bbb CP}_2$ such that $\Phi^*[g]$
is the conformal class of the Fubini-Study metric. 
\end{main} 

The basic estimate, however, has other interesting 
consequences, notably for the theory of almost-K\"ahler
manifolds.

\section{The Basic Estimate}
\label{key}

Let $M^4$ be a compact oriented 4-manifold with $b^+=1$. 
The open cone  
$$\{ [\omega]\in H^2(M, {\Bbb R})| [\omega]\cdot [\omega ] > 0\}$$ 
then consists of two connected components,
called {\em nappes}.  Given such a nappe ${\cal C}^+$,
and a Riemannian metric    $g$ on $M$,
let $\omega$ be a $g$-harmonic 2-form
such that $[\omega ]\in {\cal C}^+$; 
of course, 
such a form   always exists, and it is 
uniquely determined once we also 
specify the positive real number 
$[\omega ]^2$. 
  If $c$ is the spin$^c$ structure
induced by some almost-complex structure $J$ on $M$, then,
relative to any metric $g$, we can  
define the (perturbed) Seiberg-Witten invariant  
$p_c(M , {\cal C}^+)$ as the number of solutions,
modulo gauge and  counted with orientations, of 
the perturbed Seiberg-Witten equations \cite{KM,witten}
\begin{eqnarray} D_{A}\Phi &=&0 \label{dir} \\
 iF^+_A+\sigma (\Phi ) &=& \varepsilon , \label{sd}\end{eqnarray}
where $\varepsilon$ is a generic  self-dual 2-form 
with $\int_M \varepsilon \wedge \omega  > 2\pi c_1\cdot [\omega ]$.
Here $\omega$ is a $g$-self-dual 2-form with 
$[\omega ]\in {\cal C}^+$. More generally, if
 $c$ is  a spin$^c$ structure 
for which  $k=[c_1^2 - (2\chi + 3\tau)(M)]/4$ is non-negative 
and even, one can 
define  \cite{liu,t2}  the perturbed Seiberg-Witten invariant 
$p_c(M , {\cal C}^+)$ to be $\int_X \alpha^{k/2}$, where
$X$ is the moduli space of solutions of 
(\ref{dir}--\ref{sd}) for 
a generic $\varepsilon$
with $\int \varepsilon \wedge \omega  > 2\pi c_1\cdot [\omega ]$,
and $\alpha \in H^2(X)$ is the first Chern class of
the based moduli space, considered as an $S^1$-bundle over $X$.
 For our purposes, the key point 
is simply that when this invariant is non-zero,
  the equations 
\begin{eqnarray} D_{A}\Phi &=&0 \label{specdir}\\
 -iF^+_A&=& \sigma (\Phi ) - t\omega , \label{specsd} \end{eqnarray}
have a solution with $\Phi\not\equiv 0$
 for any $t\gg 0$.

 \begin{thm} \label{gest}
Let $M^4$ be a smooth compact oriented 4-manifold with $b^+=1$
for which there is a nappe ${\cal C}^+$ 
and a spin$^c$ structure $c$ such that 
$p_c(M,{\cal C}^+)\neq 0$. Let $g$ be any 
Riemannian metric on $M$, and let $\omega$ be a
$g$-self-dual harmonic 2-form with $[\omega]\in {\cal C}^+\subset
H^2(M,{\Bbb R})$. Then the scalar curvature $s$ of $g$
satisfies
$$\int s \frac{|\omega |}{\sqrt{2}} d\mu \leq 4\pi c_1\cdot [\omega
].$$
Here $d\mu$  and  $|\cdot |$ are respectively the 
volume form and 
  point-wise norm     determined by the metric $g$,
while $c_1=c_1(V_+)$ is the first Chern class of 
the spin$^c$ structure. 
\end{thm}
\begin{proof}
The Dirac equation $D_A\Phi =0$ implies the   Weitzenb\"ock formula
$$0=  \langle \Phi, \Delta \Phi \rangle +
 \frac{s}{4} |\Phi |^2 + 2 \langle -iF^+_A , \sigma (\Phi )\rangle
,$$
where the natural real-quadratic map
$\sigma : V_+\to \wedge^+$ satisfies $|\sigma (\Phi)|= 2^{-3/2}|\Phi |^2$.
For a solution of the perturbed Seberg-Witten equations
(\ref{specdir}--\ref{specsd}), 
 we also have $-iF^+_A = \sigma (\Phi ) - t\omega$, so it follows that
\bea 
0 &=& \int_M \left[ 4\langle \Phi, \Delta \Phi \rangle +
s |\Phi |^2 + 8 \langle  \sigma (\Phi ) - t\omega , \sigma (\Phi
)\rangle \right] d\mu
\\&=& \int_M \left[ |2\nabla_A \Phi |^2 +
s |\Phi |^2 + |\Phi |^4 - 8t \langle  \omega , \sigma (\Phi )\rangle
\right] d\mu
\\&\geq & \int_M \left[ 
s |\Phi |^2 + |\Phi |^4 - 8t | \omega | ~  |  \sigma (\Phi) | \right]
d\mu
\\&=& \int_M \left[  
s |\Phi |^2 + |\Phi |^4 - 2\sqrt{2} t | \omega | ~  |   \Phi  |^2
\right] d\mu .
\eea 
Hence 
$$ \int_M (2\sqrt{2} t | \omega |  - s) |\Phi |^2 d\mu  \geq  \int_M
|\Phi |^4 d\mu  $$
and so    the Cauchy-Schwartz inequality yields 
\bea
8t^2 [\omega ]^2 - 4\sqrt{2}t \int_M s|\omega | d\mu + \int_M s^2d\mu
&=&\int (2\sqrt{2} t  | \omega | - s)^2 d\mu 
\\&\geq & 
\frac{\left[\int    |\Phi |^4 d\mu\right]^2}{\int    |\Phi |^4 d\mu}
\\&=& 
\int    |\Phi |^4 d\mu\\  &=&8\int    |\sigma(\Phi) |^2 d\mu
\\&=& 8\int   |-iF^+_A+t\omega |^2 d\mu
\\&=& 8\int   \left( t^2|\omega |^2   -2t\langle iF^+_A,\omega
\rangle 
+  |iF^+_A|^2\right) d\mu
\\&= &  8\left(t^2  [\omega ]^2 - 2 t (2\pi c_1) \cdot [\omega ]
+  \int |iF^+_A|^2d\mu  \right)
\\&\geq& 8t^2 [\omega ]^2 - 32\pi t c_1 \cdot [\omega ] + 32\pi^2
c_1^2 .
\eea
for all $t \gg 0$.   Hence
$$32\pi t c_1 \cdot [\omega ] +\int_M s^2d\mu \geq 
4\sqrt{2}t \int_M s|\omega | d\mu + 32\pi^2c_1^2 $$
for all $t\gg 0$. Dividing by $8t$ and taking the limit as
$t\to + \infty$, we therefore have
$4\pi c_1\cdot [\omega ] \geq \int_M s \frac{|\omega
|}{\sqrt{2}}d\mu$,
as claimed. 
\end{proof}

\begin{thm}\label{eq} 
Equality holds in Theorem  \ref{gest} 
iff $g$ is  K\"ahler  
with respect to some $c$-compatible complex structure  $J$ and  
 $\omega$ is a constant positive 
multiple of the K\"ahler form of $(M,g,J)$. 
\end{thm}
\begin{proof}
The Weitzenb\"ock formula 
for  (\ref{specdir}--\ref{specsd}) tells us that
\bea
0&=&2 \Delta |\Phi |^2 +
 |2\nabla \Phi |^2 +  s  |\Phi |^2 +  |\Phi |^4-  8t \langle \omega ,
\sigma (\Phi )\rangle\\
&\geq &   2 \Delta |\Phi |^2 +
s  |\Phi |^2 +  |\Phi |^4 - 2\sqrt{2} t |\omega | |\Phi |^2,
\eea
so that 
$$|\Phi |^2 \leq 2\sqrt{2} t |\omega | - s$$
at the maximum of $|\Phi |$. The $C^0$ norm of the 
twisted spinor field $\Psi_t = \Phi /\sqrt{t}$
 is therefore uniformly bounded as $t\to \infty$: 
$$|\Psi_t |^2 \leq 2\sqrt{2} \max |\omega  | -\min (\min s ,0)
~\forall t > 1.$$

Now the   proof of Theorem \ref{gest}     contains the
inequality
$$8t^2 [\omega ]^2 - 4\sqrt{2}t \int_M s|\omega | d\mu + \int_M
s^2d\mu 
\geq 
 8\left(t^2  [\omega ]^2 - 2 t (2\pi c_1) \cdot [\omega ]
+ \int_M |iF^+_{A_t}|^2d\mu\right) 
$$
and it therefore   follows that 
$$
\int_M s\frac{|\omega |}{\sqrt{2}} d\mu =4\pi  c_1 \cdot [\omega ] ~
\Longrightarrow ~ \|iF^+_A\|_{L^2}^2 \leq \frac{1}{8} \int_Ms^2d\mu .
$$
Since equation (\ref{specsd}) stipulates that 
$$\sigma (\Psi_t) -\omega = 
\frac{\sigma (\Phi ) - t\omega}{t}   
= - \frac{iF_A^+}{t}, $$
equality in Theorem \ref{gest} implies
 that $\sigma (\Psi_t) \to \omega$ in $L^2$.

On the other hand, the proof of Theorem \ref{gest} 
neglected the  $|\nabla \Phi |$ term. 
Leaving it in until the    Schwartz inequality step 
yields 
$$\int (2\sqrt{2} t  | \omega | - s)^2 d\mu 
 \geq  
\frac{\left[\int |2\nabla_A \Phi |^2d\mu +
\int    |\Phi |^4 d\mu\right]^2}{\int    |\Phi |^4 d\mu}
 \geq  
\int    |\Phi |^4 d\mu + 8\int |\nabla_A \Phi |^2 d\mu , $$
and one therefore concludes that 
$$t\left[4\pi  c_1 \cdot [\omega ] - \int_M s\frac{|\omega
|}{\sqrt{2}} d\mu
\right]
+\frac{1}{8}\left[\int_M s^2d\mu  
 - 32\pi^2c_1^2 \right]\geq \int_M |\nabla_A \Phi |^2 d\mu ,$$
which we may rewrite as 
$$\left[4\pi  c_1 \cdot [\omega ] - \int_M s\frac{|\omega |}{\sqrt{2}}
d\mu
\right]
+\frac{1}{8t}\left[\int_M s^2d\mu  
 - 32\pi^2c_1^2 \right]\geq \int_M |\nabla_{A_t} \Psi_t |^2 d\mu.$$
 It follows that 
$$ \int_M s\frac{|\omega |}{\sqrt{2}} d\mu =4\pi  c_1 \cdot [\omega ]
\Longrightarrow ~ \lim_{t\to \infty}\|\nabla_{A_t}  \Psi_t\|_{L^2}
=0.$$
Since $\sigma (\Psi_t)$ is the contraction of a parallel
field with $\Psi_t\otimes\bar{\Psi}_t+\bar{\Psi}_t \otimes\Psi_t$, and
because 
we   have a uniform $C^0$ bound for $|\Psi_t|$, application of
the Leibnitz rule   yields 
$$\lim_{t\to \infty}\|\nabla \sigma(\Psi_t)\|_{L^2} =0$$
whenever equality holds in Theorem \ref{gest}. 
Here the connection $\nabla$ on self-dual 2-forms is 
  $t$-independent, and is actually the one 
induced by  the (torsion-free) Levi-Civit\`a connection of $g$. 
Since   $\omega$ is harmonic, it follows that 
$$\lim_{t\to \infty}\| (d^+)^* ( \sigma(\Psi_t)-\omega )\|_{L^2} =
\lim_{t\to \infty}\| \nabla\cdot \sigma(\Psi_t) \|_{L^2} =0,$$
where $(d^+)^*= -\star d$
is the adjoint of $d^+ : \Gamma\wedge^1\to  \Gamma\wedge^+$ with
respect 
to $g$. 
Because the complex 
$$0\to \Gamma \wedge^0 \stackrel{d}{\longrightarrow} \Gamma \wedge^1
 \stackrel{d^+}{\longrightarrow}  \Gamma \wedge^+\to 0$$
is elliptic, the G{\aa}rding inequality 
$$\|  \phi\|_{L^2_1} \leq C\left(\|  (d^+)^* \phi\|_{L^2}+ \|
\phi\|_{L^2}
\right)  ~~ \forall \phi \in \Gamma \wedge^+ , $$
applied to $\phi = \sigma(\Psi_t)-\omega $, 
then tells us that $\sigma (\Psi_t) \to \omega$ in the Sobolev space
$L^2_1$. It follows that 
$$\|\nabla \omega \|_{L^2} = \lim_{t\to \infty} \|\nabla \sigma
(\Psi_t) \|_{L^2} =0, $$
and hence that $\nabla\omega =0$. Thus equality in
Theorem \ref{gest} implies that $g$ is K\"ahler, as claimed. 
\end{proof}

By the same reasoning, we also get 

 \begin{thm}
Let $M^4$ be a smooth compact oriented 4-manifold with $b^+ >1$
for which there is   a spin$^c$ structure $c$ such that 
the Seiberg-Witten invariant 
$n_c(M)$ is non-zero. Let $g$ be any 
Riemannian metric on $M$, and let $\omega$ be any
$g$-self-dual harmonic 2-form. Then the scalar curvature $s$ of $g$
satisfies
$$\int s \frac{|\omega |}{\sqrt{2}} d\mu \leq 4\pi c_1\cdot [\omega
].$$
If $[\omega ]\neq 0$, moreover, equality is achieved  iff $g$ is
K\"ahler  
with respect to some $c$-complex structure  $J$ and  
 $\omega$ is a constant positive 
multiple of the K\"ahler form of $(M,g,J)$. 
\end{thm}

Recall \cite{goldberg} that an {\em almost-K\"ahler manifold} is
a triple
$(M,g,J)$, where $J$ is an almost complex structure on  $M$ and 
$g$ is a $J$-invariant Riemannian metric, such that 
the 2-form $\omega (\cdot , \cdot ) = g(J\cdot, \cdot )$
is closed. The symplectic form $\omega$ is then 
called the {\em  almost-K\"ahler form} of
 $(M,g,J)$. 
 
\begin{cor} \label{ak}
Let $(M^4,g,J )$ be an almost-K\"ahler manifold. 
Then the scalar curvature $s$ of $g$ satisfies 
$$\int_M s ~d\mu \leq 4\pi c_1\cdot [\omega ],$$
where $c_1$ is the first Chern class of  $(TM, J)$
and $\omega$ is the almost-K\"ahler form. 
Moreover, equality is achieved iff  $(M,g,J)$
is K\"ahler. 
\end{cor}
\begin{proof}
On any almost-K\"ahler manifold $(M^4,g,J)$, 
the almost-K\"ahler form is a harmonic self-dual 
form with $|\omega |\equiv \sqrt{2}$. Since Taubes 
\cite{taubes} has shown 
that the perturbed Seiberg-Witten invariant of any symplectic
4-manifold is non-zero, the result therefore  follows from the
preceding theorems. 
\end{proof}

Note that  the  associated variational problem
was first studied in \cite{blair}.

\section{Genericity and Self-Duality}

In light of Corollary \ref{ak}, it is easy to 
see that Theorem \ref{gest} gives rise 
to an estimate of Yamabe invaraints whenever
the harmonic 2-form has empty zero locus.
To generalize this, we will need to 
understand the manner in which harmonic forms
vanish for generic metrics. 
  The needed information is provided by the
following  result, the essence of which is contained 
 in unpublished work  of Taubes \cite[Lemma B.3]{t2}.

\begin{prop}[Taubes]  \label{tap}
Let $M$ be a smooth oriented 4-manifold
with $b_+\geq 1$, and let $\cal M$ denote the space of 
$C^{k,\alpha}$ Riemannian 
metrics $g$ on $M$ for some  $k\geq 2$, $\alpha \in (0,1)$.  
Choose a codimension-$(b_+-1)$ subspace 
$V\subset H^2(M, {\Bbb R})$ on which the 
intersection form has Lorentzian signature, 
and let $V_1$ be a connected component of
the hyperboloid
$\{ [\omega ]\in V~|~ [\omega ]^2 =1 \}$. 
For each $g \in {\cal M}$, let $\omega_g$ be 
the unique self-dual $g$-harmonic form with
$[\omega_g ]\in V_1$. 
Then  
$$\{ ~  g ~|~   
\omega_g  \mbox{ is transverse to the zero section of }
 \wedge^+_g \} \subset {\cal M}$$
 is a subset of the second Baire category. 
\end{prop}
\begin{proof} 
Let  $\bigwedge^+\to M\times {\cal M}$ 
be the  
$C^k$ rank-$3$ vector bundle  whose restriction 
to any slice  $M\times \{  g\}$ is the bundle
$\wedge^+_g\to M$
 of  self-dual 2-forms determined by
 the metric $g$. Let $\Omega$ denote the 
section of $\bigwedge^+$
defined by  
$$\Omega |_{M\times \{  g\}}= \omega_g.$$ 
By a   straight-forward application of the 
Banach-space inverse function theorem
and the interior Schauder estimates \cite{DN}, 
the map ${\cal M}\ni g\mapsto \omega_g \in  C^{k,\alpha}(\wedge^2)$
  can  be seen to be smooth,
 and this implies that $\Omega$ is  a $C^k$
section of $\bigwedge^+$. 

Now $\Omega$  turns out to  actually be transverse to the
zero section of $\bigwedge^+$. 
To see this, let $g\in {\cal M}$ and 
observe that for any other $\tilde{g}\in {\cal M}$,
the bundle
$\wedge^+_{\tilde{g}}\subset \wedge^2$
may be thought of as the graph of a homomorphism
$h:\wedge^+_g\to \wedge^-_g$, and so in particular
may be identified with $\wedge^+_g$ via 
the projection $\wedge^2\to \wedge^+_g$. 
If $\tilde{g}$ is conformally rescaled so as to have
the same volume form of $g$, moreover, 
and is we identify ${\cal H}om (\wedge^+, \wedge^-)$
with $\odot_0^2T^*M$ in the usual way, then one has
$\tilde{g} = g + h + O(|h|^2)$. 
Using  these identifications, 
one can then check that 
$$\left.
\frac{d}{dt}\omega_{g+th}\right|_{t=0}=  -2  G [ d^+\delta h(\omega )
]  $$
by simply using only the facts that $\omega_{\tilde{g}}$ is 
closed and depends differentiably on $\tilde{g}$;
here $G$ is the Green's operator of the Hodge Laplacian 
$d\delta + \delta d$
of $g$ on
2-forms.  If $\omega$ vanishes at $p\in M$, it therefore
suffices to show that $h\in \Gamma {\cal H}om (\wedge^+, \wedge^-)=
\Gamma \odot_0^2T^*M$ may be chosen so that 
$G[ d^+\delta h(\omega ) ]|_p$ is any desired element   
$\phi$ of 
$\wedge^+_g|_p$.  
If $\hat{\phi}$ is the  
delta-function-like 2-current  
$\hat{\phi}(\psi )=\langle \psi|_p , \phi \rangle_g$
corresponding to a non-zero  $\phi\in \wedge^+_g|_p$,
however, 
this is equivalent to the assertion that the current
$d^-\delta G \hat{\phi}$ is not identically zero on 
the open  subset of $M$ where $\omega_g\neq 0$.
Now the latter open subset is also dense in $M$ by
unique continuation. 
Since $G \hat{\phi}$ is $C^2$ on $M-p$
by elliptic regularity, it therefore suffices to show
that $d \delta G \hat{\phi}$ is not identically
zero on $M-p$.   But in geodesic normal coordinates
$x^j$
about $p$, one has an expansion \cite{combet} 
$$G(\hat{\phi})=
\frac{1}{4\pi^2} 
\frac{\phi}{|x|^2} + O(\log |x| )$$
where the Hessian of the $O(\log |x| )$ term is
$O(1/|x|^2)$; here $\phi$ has been extended as 
tensor field such  that $\nabla \phi |_p =0$. 
Assuming without loss of generality that 
$|\phi |  = \sqrt{2}$ and rotating our coordinates
as necessary, we may arrange, after 
transplanting the standard complex structure $J$
and complex codinates $(z^1, z^2)$ on 
on ${\Bbb R}^4={\Bbb C}^2$ to our coordinate
chart,  that $\phi  =  g(J\cdot, \cdot )$,
so that 
\bea d\delta G(\hat{\phi})&=&
-\frac{1}{4\pi^2} 
dJd (\frac{1}{r^2})  + O( \frac{1}{r^2})\\
&=&\frac{i(|z^2|^2-|z^1|^2) (dz^1\wedge d\bar{z}^1
-dz^2\wedge d\bar{z}^2)+ 4 \Im m  
	(\bar{z}^1{z}^2 dz^1\wedge d\bar{z}^2)}{2\pi^2r^6}
+ O( \frac{1}{r^2})
\eea
is indeed non-zero, and 
   $\Omega$ is transverse to the zero section, as claimed.

Consequently, the zero locus  $Z\subset M\times {\cal M}_{[\omega ]}$
  of $\Omega$ is a $C^k$ Banach submanifold of  
 codimension $3$, and  the induced projection 
$\wp : Z\to {\cal M}_{[\omega ]}$ is    Fredholm,
of  index $4-3=1$. Since $\wp$ is a $C^k$ map,
and $k > \max (\mbox{index}(\wp) , 0)$, 
 the Smale-Sard Theorem \cite{sm} tells us  that the  
set of regular values of $\wp$ is a set of the second Baire
category.
 But   this  set of regular values is 
precisely the set of $C^{k,\alpha}$ metrics 
for which $\omega_g$ is transverse to the zero section.
\end{proof}

This immediately implies the technical result
 we'll actually use:
 
\begin{cor} \label{tac} 
Let $g$ be a  $C^2$ Riemannian metric
on a compact oriented  4-manifold $M$,
 and 
let $\omega \not\equiv 0$ be a self-dual harmonic form 
on $(M,g)$. Then for any integer $r \geq 2$ there is a sequence 
$\{ g_j\}$ of $C^r$ metrics
 on $M$ and a sequence $\{\omega_j \}$ of closed 2-forms  such that 
\begin{description}
\item{(a)} $\lim_{j\to \infty} g_j =g$ in the $C^2$ topology;
\item{(b)} $\omega_j$ is self-dual with respect to $g_j$;
\item{(c)}   $\omega_j$
is transverse to the zero section of $\wedge^+_{g_j}\to M$; and
\item{(d)} $\lim_{j\to \infty} \omega_j =\omega$ in the $C^1$
topology.
\end{description}
\end{cor}
\begin{proof}
Let
$V\subset H^2(M, {\Bbb R})$ be   
spanned by $[\omega ]$ and the cohomological image
of the anti-self-dual 
$g$-harmonic forms; and let $V_1$  be the sheet  
of the unit hyperboloid in $V$ which contains $\lambda [\omega ]$ for
some
$\lambda  > 0$. For each $C^2$ metric $\hat{g}$ we then have a 
unique self-dual harmonic 2-form $\omega_{\hat{g}}$ with
$[\omega_{\hat{g}}
 ]\in V_1$. 
By Proposition \ref{tap}, $\omega_{\hat{g}}$ is transverse to the 
zero section for  a dense set of $C^{2,\alpha}$ metrics $\hat{g}$, 
and this set is automatically
also dense in the space of $C^2$ metrics. Let $g_j$ be a sequence of
such metrics converging to $g$ in the $C^2$ topology,
and set $\omega_j=\omega_{g_j}/\lambda$. 
Now,
for any 
$\alpha \in (0,1)$,
 $\hat{g}\mapsto \omega_{\hat{g}}$ is a smooth map from 
$C^{1,\alpha}$ metrics to $C^{1,\alpha}$ forms,
  so $\omega_j \to \omega$ in the 
$C^{1,\alpha}$ topology, and hence in the $C^1$ topology,
as claimed.  
\end{proof}

\section{Yamabe Invariants}
\label{rep}

Let $M$ be a smooth compact oriented 4-manifold, and let 
$[g]$ be a conformal class of metrics on $M$. The Yamabe
constant of the conformal class is then  defined by
$$Y_{[g]} = \inf_{g\in [g]}\frac{\int_M s_g ~ d\mu_g}{\sqrt{\int_M
d\mu_g}}.$$
Since Hodge star operator $\star : \wedge^2\to \wedge^2$
is conformally invariant, the corresponding eigenspace decomposition 
$$\wedge^2= \wedge^+\oplus \wedge^-$$
of the 2-forms into self-dual and anti-self-dual
parts is independent of a choice of $g\in [g]$.

\begin{thm}\label{yam}
Let $(M,[g])$ be an oriented  conformal Riemannian 4-manifold,
and let $\omega\not\equiv 0$ be a closed  2-form which is self-dual  
with respect to $[g]$.  
 Suppose that  $b^+(M) = 1$ and that the perturbed Seiberg-Witten
invariant $p_c(M,{\cal C}^+)$ is non-zero for some spin$^c$ structure
$c$,
where ${\cal C}^+\subset H^2(M,{\Bbb R})$ is the nappe containing
$[\omega ]$.  
Then the Yamabe constant of
$[g]$ satisfies 
$$Y_{[g]}\leq \frac{4 \pi c_1\cdot[\omega]}{\sqrt{[\omega]^2/2}},$$
where $c_1$ is the first Chern class of $c$. 
\end{thm}
\begin{proof}
We may assume that the graph of the $[g]$-harmonic 2-form  
$\omega$ is transverse to the zero section of $\wedge^+$. 
Indeed, Corollary \ref{tac} tells us that the set of 
conformal classes with this property
is dense, and both sides of the inequality are
 continuous functions 
with respect to the $C^2$-topology \cite{ber}.

Let $ \varphi: {\Bbb R}\to {\Bbb R}$ be a smooth   function
  with 
\begin{itemize}
\item $\varphi(x)= \frac{1}{2}$ when $x\leq \frac{1}{4}$;
\item $\varphi (x)=x$ when $x\geq 1$; and 
\item $0\leq \varphi^{\prime} (x) \leq 1$ for all $x$.
\end{itemize} 
Let $g$ be any representative of $[g]$, 
and for each $\delta > 0$ define a smooth
positive function by  
$v=v_{\delta}=\delta \varphi ( |\omega |_g/\delta ) \geq |\omega
|_g$.
Let 
$$U_\delta =\{ p\in M ~|~ |\omega |_g \leq \delta \}. $$ 
Notice that 
\begin{itemize}
\item $|dv|_g < C$;
\item $\frac{\delta}{2} \leq v \leq  \delta$ on
$U_\delta$; 
\item  $v\equiv  \frac{\delta}{2}$   on $U_{\delta/4}$;
and 
\item $v \equiv |\omega |_g$ on $M-U_\delta$, 
\end{itemize}
where  $C$ is some constant independent of $\delta$.
Moreover, 
 there is a number $\delta_0 \in (0,1)$ such that
for all $\delta \in (0, \delta_0 )$  
\begin{itemize}
\item   $U_{\delta}$   is a smooth manifold with boundary; and
\item   $\mbox{Vol}_g(U_{\delta})< C\delta^3$.   
\end{itemize}
 Indeed,  
we may take $C= \max (\max |\nabla\omega |_g , 2 \pi \ell )$, 
where $\ell$ is the  total length  of the compact curve
 $\omega =0$, measured with respect to $g$.  Notice
that we are now exploiting the hypothesis that $[g]$
is generic. 

For any $\delta \in (0, \delta_0)$,   consider the smooth 
conformal metric $\hat{g}=vg\in [g]$. 
One has 
$$|\omega |_{\hat{g}} = v^{-1} |\omega |_{\hat{g}}\leq
1$$
and $|\omega |_{\hat{g}}\equiv 1$ on $M-U_{\delta}$. 
Moreover, the continuous function  $\psi = |\omega|_{\hat{g}}: M\to
[0,1]$
is smooth on $M-U_{\delta/4}$, with $|d\psi|_g \leq 20 C\delta^{-2}$
on this set, 
so  
\bea 
\left|  \int_Ms_{\hat{g}}d\mu_{\hat{g}} 
-\int_M s_{\hat{g}} |\omega|_{\hat{g}}d\mu_{\hat{g}} \right| &=& 
\left| \int_{U_{\delta}}
 s_{\hat{g}} (1-\psi)d\mu_{\hat{g}} \right|
\\&= & \left| \int_{U_{\delta}}  (1-\psi)
\left[3\Delta v + s_g v + 
\frac{3}{2} \frac{|dv|^2_g}{v}\right]
d\mu_g   \right|
\\&\leq  &3 \left| \int_{U_{\delta}}  (1-\psi)\Delta v
d\mu_g\right|
\\&&~~+  \left| \int_{U_{\delta}}  (1-\psi)
s_g v 
d\mu_g   \right|
+  \frac{3}{2}\left| \int_{U_{\delta}}  (1-\psi)
 \frac{|dv|^2_g}{v}
d\mu_g   \right|
\\&\leq  &3 \left| \int_{U_{\delta}-U_{\delta /4}}  
(1-\psi)\Delta v d\mu_g\right|
\\&&~~+    \int_{U_{\delta}}   
|s_g |v 
d\mu_g  
+  \frac{3}{2}  \int_{U_{\delta}}   
 \frac{|dv|_g^2}{v}
d\mu_g    
\\&\leq  &3 \left| \int_{\partial U_{\delta}\cup \partial
\bar{U}_{\delta /4}} 
(1-\psi)\star dv +
  \int_{U_{\delta}-U_{\delta /4}}  
\langle d\psi, dv \rangle_g
d\mu_g    
\right|
 \\&&~~
+      (\max |s_g |) C\delta^4
+   C^2\frac{3}{\delta} C\delta^3
\\&\leq  & 3C\frac{20 C}{\delta^2} C\delta^3
+     C (\max |s_g |) \delta^4
+  3 C\delta^2 \leq \hat{C}\delta
\eea 
because $\psi \equiv 1$ on $\partial U_\delta$ and
$dv \equiv 0$ on $\partial U_{\delta/4}$.
On the other hand,
$$\int_Md\mu_{\hat{g}}\geq  \int_M|\omega |_{\hat{g}}^2d\mu_{\hat{g}}
=\int_M \omega\wedge\omega =[\omega]^2.$$
Thus,
given $\epsilon > 0$,  consider the 
metric $\hat{g}$ corresponding to  
 $\delta= \epsilon\hat{C}^{-1}\sqrt{[\omega]^2}$.
By Theorem \ref{gest}, this metric satisfies 
$$\frac{\int_Ms_{\hat{g}}d\mu_{\hat{g}}}{\sqrt{\int_Md\mu_{\hat{g}}}}
\leq   \frac{ \int_Ms_{\hat{g}}
|\omega |_{\hat{g}}d\mu_{\hat{g}}}{\sqrt{[\omega]^2}}  + \epsilon
\leq  \frac{  4\sqrt{2}\pi c_1\cdot[\omega ]}{\sqrt{[\omega]^2}} +
\epsilon.$$
Hence 
$$Y_{[g]} = \inf_{\hat{g}\in
[g]}\frac{\int_Ms_{\hat{g}}d\mu_{\hat{g}}}{\sqrt{\int_Md\mu_{\hat{g}}}
}
\leq  \frac{  4 \pi c_1\cdot[\omega ]}{\sqrt{[\omega]^2/2}}  ,$$
as claimed.
\end{proof}
 
\begin{thm}
Equality is achieved in  Theorem \ref{yam} iff 
there is a Yamabe minimizer $g\in [g]$ which
is K\"ahler, with K\"ahler form $\omega$. 
\end{thm}
\begin{proof}
We may assume that $c_1\cdot [\omega ] > 0$, since
otherwise the result follows from the
 theory of the  unperturbed
Seiberg-Witten equations \cite{leb2}.

Let us begin by showing that equality can only be 
achieved if the harmonic 
form  $\omega$ is nowhere zero. 
Since $Y_{[g]} >0$,  
we may therefore choose  a metric
$g$ of positive scalar curvature in $[g]$.
By Corollary \ref{tac}, we can find 
  be a sequence $g_j$ of positive-scalar-curvature 
metrics  converging  to $g$ in the $C^2$ topology
such   that 
the corresponding harmonic forms $\omega_j$,
normalized so that $[\omega_j]^2=[\omega ]^2$, 
 are transverse to 
the zero section and 
converge to $\omega$ in the 
$C^1$ topology. For each $j$ there is an 
$\epsilon_j > 0$ such that 
$$
\frac{\int [s_g u^2  + 6 |du  |^2_g ]d\mu_{g}}{\sqrt{\int u^4 d\mu_{
g}}}
\leq (1+ \epsilon_j)
\frac{\int [s_{g_j} u^2  + 6 |du  |^2_{g_j} ]d\mu_{g_j}}{\sqrt{\int
u^4 d\mu_{ g_j}}}
$$
for every $u\in L^2_1$, and moreover one can take $\epsilon_j\to 0$
as $j\to \infty$. Now let $u_j=(\delta_j \varphi
(|\omega_j|_{g_j}/\delta_j )^{1/2}$,
where $\varphi$ is as in the proof of Theorem \ref{yam} and 
$\delta_j$ is chosen so small that
$$\frac{\int [s_{g_j} u^2  + 6 |du  |^2_{g_j} ]d\mu_{g_j}}{\sqrt{\int
u^4 d\mu_{ g_j}}}
\leq  \frac{4\pi c_1\cdot [\omega_j]}{\sqrt{[\omega ]^2/2}} +
\epsilon_j,$$
for each $j$, and 
is also chosen in such a manner  that $\lim_{j\to \infty} \delta_j =0$.  
Since $[\omega_j]\to [\omega ]$, it   follows that 
$$
\frac{\int [s_g u_j^2  + 6 |du_j  |^2_g ]d\mu_{g}}{\sqrt{\int u_j^4
d\mu_{ g}}}
\longrightarrow
\frac{4\pi c_1\cdot [\omega_j]}{\sqrt{[\omega ]^2/2}}=Y_{[g]}
$$
and the $u_j$ are in particular bounded in $L^2_1$ because 
  $\int u^4_j d\mu  \to [\omega]^2$. By passing to a 
subsequence, we may therefore assume that the $u_j$ are
weakly convergent to some ${u}\in L^2_1$. 
But this implies that ${u}$ is the unique weak limit of the 
$u_j$ in $L^2$. However,  
 by construction, $u_j\to \sqrt{|\omega |_g}$ in $C^0$,
and hence in $L^2$. We therefore have  $u= \sqrt{|\omega |_g}$,
so that  $ \sqrt{|\omega |_g}\in L^2_1$. 
Now since $u$ is the weak limit of the $u_j$,
$$\| u \|_{L^2_1} \leq \lim_{j\to \infty}  \| u_j \|_{L^2_1},$$
whereas we also know that 
$$\|u\|_{L^4}^4=\|\sqrt{|\omega |}\|_{L^4}^4 = [\omega ]^2 =
\lim_{j\to \infty} \|u_j\|_{L^4}^4$$
by construction. Hence
$$
\frac{\int [s_g u^2  + 6 |du  |^2_g ]d\mu_{g}}{\sqrt{\int u^4 d\mu_{
g}}}
\leq \lim_{j\to \infty} 
\frac{\int [s_g u_j^2  + 6 |du_j  |^2_g ]d\mu_{g}}{\sqrt{\int u_j^4
d\mu_{ g}}}
=Y_{[g]} = \inf_{\hat{u}\in L^2_1-0}
\frac{\int [s_g \hat{u}^2  + 6 |d\hat{u}  |^2_g ]d\mu_{g}}{\sqrt{\int
\hat{u}^4 d\mu_{ g}}}$$
so that $u\in L^2_1$ is a minimizer of the Yamabe functional,
and  is automatically a weak solution of the associated
Euler-Lagrange equation
\be \label{el} 
6\Delta u + s u = s_0 u^3 ,
\ee
where $s_0=4\sqrt{2}\pi\frac{c_1\cdot [\omega ]}{[\omega ]^2}$. 
Since    $u=\sqrt{|\omega |}\in C^{0,1/2}$, 
elliptic regularity  \cite{aubin}, applied to (\ref{el}), 
implies  that   
$u\in C^{2,1/2}$, and an application \cite[Prop. 3.75]{aubin}
of the  maximum principle to (\ref{el})
therefore shows that $u\geq 0$ must either 
be strictly positive or vanish identically. 
Hence $|\omega | = u^2\neq 0$, and 
$\omega \neq 0$, as claimed.

We may therefore define an almost K\"ahler metric $\hat{g}\in [g]$
by $\hat{g}= \sqrt{2} |\omega|_g  g$. Since $\int d\mu_{\hat{g}}
=[\omega ]^2/2$, Corollary \ref{ak} tells us that 
$$\frac{\int s_{\hat{g}}~ d\mu_{\hat{g}}}{\sqrt{\int d\mu_{\hat{g}}}}
\leq 
\frac{4\pi 
c_1\cdot [\omega ]}{\sqrt{[\omega]^2/2}},$$
with equality  iff $g$ is K\"ahler.  Since
the right-hand side is  $Y_{[g]}$ by hypothesis, 
however, we also have a tautological inequality 
in the other direction, and the two sides must be equal.
Thus $\hat{g}$ is both a Yamabe minimizer and a  K\"ahler 
metric, as desired.
\end{proof}

Now recall that the {\rm Yamabe invariant} (sometimes called the
 sigma constant)  of a 
smooth $n$-manifold $M$ is defined to be 
$$Y(M)= \sup_{[g]}
\inf_{g\in [g]}\frac{\int_Ms_gd\mu_g}{\left( \int_Md\mu_g
\right)^{\frac{n-2}{n}}}
= \sup_{[g]} Y_{[g]}$$
where $[g]$ is allowed to vary over the  space of conformal class  of
smooth  Riemannian
metrics on $M$. 

\begin{thm}[Main Theorem] \label{cp2}
The Yamabe invariant of the complex projective plane is given by
$$Y({\Bbb CP}_2) = 12\sqrt{2}\pi .$$
Moreover, a conformal class $[g]$ on ${\Bbb CP}_2$ satisfies
$Y_{[g]} = Y({\Bbb CP}_2)$ iff there is a diffeomorphism 
$\Phi : {\Bbb CP}_2\to {\Bbb CP}_2$ such that $\Phi^*[g]$
is the conformal class of the Fubini-Study metric. 
\end{thm} 
\begin{proof}
Let $c$ be the spin$^c$ structure induced by 
the usual complex structure on ${\Bbb CP}_2$,
and let ${\cal C}^+$ denote the nappe containing
$c_1$. Then  the perturbed Seiberg-Witten
invariant $p_c({\Bbb CP}_2, {\cal C}^+)\neq 0$,
as may be seen in a variety of ways \cite{KM,leb1,liu,taubes}.
Now since $b_-=0$,
the class $[\omega ] = c_1$ is self-dual with respect to any
conformal class $[g]$, and we therefore always have 
$$Y_{[g]}\leq 4\pi \frac{c_1\cdot c_1}{\sqrt{c_1^2/2}}=
12\sqrt{2}\pi,$$
 with equality  iff [g] is
 represented by  a K\"ahler  metric metric
of constant positive scalar curvature. 
For such a metric, however, the Ricci form is 
harmonicm and since 
  $b_-=0$,  this implies that  the Ricci form is 
a constant multiple of the K\"ahler form, and the 
metric is therefore K\"ahler-Einstein, with positive scalar
curvature.
But it follows, for example, from the Enriques-Kodaira classification
\cite{bpv} and   Matsushima's theorem on isometry 
and automorphism groups \cite{mats} that
any such metric on the smooth 4-manifold 
 ${\Bbb CP}_2$ is  the standard Fubini-Study 
metric up to  diffeomorphisms and  rescaling. 
\end{proof}

\begin{thm} \label{cp?}
Let $M$ be a smooth compact orientable manifold
with $b_1=0$ and $b_2=1$. Then  
$$Y(M) \leq   12\sqrt{2}\pi .$$
Moreover,    $M$ admits a conformal class with
$Y_{[g]}= 12\sqrt{2}\pi$ iff $M$ is diffeomorphic to  ${\Bbb CP}_2$.
\end{thm} 
\begin{proof}
If $M$ does not admit metrics of positive scalar curvature, 
there is nothing to prove. Otherwise, let $c$ be 
a spin$^c$-structure with $c_1^2= 2\chi + 3\tau = 9$,
and let ${\cal C}^+$ be the nappe containing 
$[\omega ] = c_1$. Since there is a metric 
of positive scalar curvature on $M$, there are no 
solutions of the unperturbed Seiberg-Witten equations 
for any metric on $M$, and because $c_1\cdot [\omega ] > 0$,
the wall-crossing formula \cite{KM,liu} therefore
allows us to conclude that the perturbed invariant 
$p_c(M, {\cal C^+})$ is $\pm 1$. 
As in the last proof, we therefore have
$$Y_{[g]}\leq 4\pi \frac{c_1\cdot c_1}{\sqrt{c_1^2/2}}=
12\sqrt{2}\pi$$
for every conformal class, with equality only if [g] is
 represented by  a K\"ahler-Einstein metric
of positive scalar curvature.  
The diffeomorphism statement therefore follows 
from the classification \cite{bpv} of  complex surfaces with  $c_1 >
0$. 
\end{proof}

Exactly the same reasoning also proves the following:

\begin{prop}  
Let $(M,J)$ be a compact complex surface with $c_1 > 0$,
and let $[g]$ be any conformal class 
on $M$ for which $c_1$ is self-dual.
Then 
$$Y_{[g]}\leq   4\pi \sqrt{2c_1^2(M)},$$
with equality iff $[g]$ can be represented by 
a K\"ahler-Einstein metric adapted to some
deformation of $J$. 
\end{prop}

However, this is insufficient to  compute
the Yamabe invariant   manifolds other than
${\Bbb CP}_2$. For example, the 
conformal class
$[h]$ of Page's Einstein
metric \cite{bes} on ${\Bbb CP}_2\# \overline{\Bbb CP}_2$
has $Y_{[h]} > 4\pi \sqrt{2c_1^2({\Bbb CP}_2\# \overline{\Bbb CP}_2
)}$.
This of course  does not contradict the above result, 
since   
$c_1$ is not self-dual with respect to $[h]$. 
Indeed, since there are sequences of metrics 
on any blow-up of ${\Bbb CP}_2$ which bubble off
to metrics on ${\Bbb CP}_2$, 
it seems reasonable to conjecture
that $Y({\Bbb CP}_2\# k\overline{\Bbb CP}_2 )= 
Y({\Bbb CP}_2) = 12\sqrt{2}\pi$
for all $k$.

Let us now conclude  with a new, 
non-twistor-theoretic  proof of a   theorem 
of  Poon \cite{poon}.
Recall that an oriented conformal Riemannian manifold $(M,[g])$
is called self-dual if its Weyl curvature satisfies $W_-=0$.
For example, the conformal class of the 
 Fubini-Study metric on  ${\Bbb CP}_2$  
is self-dual.

\begin{cor}[Poon]
Let $(M,[g])$ be a self-dual 4-manifold with $b_1=0$ and $b_2 =1$.
Suppose, moreover, that $Y_{[g]} \geq 0$. Then $(M,[g])$
is conformally isometric to ${\Bbb CP}_2$  equipped with
the Fubini-Study metric. 
\end{cor}
\begin{proof}
Let $g\in [g]$ be a Yamabe minimizer.
Since $$b_+-b_-= 
\tau = \frac{1}{12\pi^2} \int_M \left[ |W_+|^2-|W_-|^2\right]d\mu \geq
0, $$
and $b_++b_-=b_2=1$, 
we must have $b_+=1$ and $b_-=0$. 
Thus
$$3= 2\chi -3\tau = \frac{1}{4\pi^2}\int_M \left[ 2|W_-|^2+
\frac{s^2}{24}
 - \frac{|\stackrel{\circ}{r}|^2}{2}\right]d\mu \leq
\frac{1}{96\pi^2}\int_M s^2d\mu$$
where $\stackrel{\circ}{r}$ denotes the trace-free Ricci tensor
and $W_-=0$. 
 Because $g$ was chosen to
be a Yamabe minimizer, and so has constant   scalar curvature, 
Theorem \ref{cp?} tells us that
$$\int_M s^2 d\mu = (Y_{[g]})^2 \leq (12\sqrt{2}\pi)^2= 3(96\pi^2),$$
so the above inequality is actually   an equality,
forcing 
$$Y_{[g]}= 12\sqrt{2}\pi,$$
and so implying that $M$ is diffeomeorphic to ${\Bbb CP}_2$. 
Moreover, Theorem \ref{cp2} tells us that the
diffeomeorphism can be chosen so that  $[g]$ is the
pull-back of the Fubini-Study conformal class. 
\end{proof}

\vfill

\noindent
{\bf Acknowledgements.} 
The author would particularly like to thank
 Cliff Taubes for his 
help in connection with Proposition \ref{tap}. 
He would also like to thank Mike Anderson and Matt 
Gursky for their helpful remarks on minimizing sequences in 
the Yamabe problem.   
 
\pagebreak

\end{document}